# Polya by Examples

KUNG-WEI YANG

E-mail: kung-wei.yang@wmich.edu

Polya's Enumeration Theorem is one of the most useful tools dealing with the enumeration of patterns that are symmetric in some ways. What follows is a procedure for obtaining the *results* of Polya's Theorem directly, bypassing the usual preliminaries ("cycle index,"...). It is for people who want to *use* (and to understand) Polya's Theorem. There is no proof here. There is no fancy example here.

"Always begin with the simplest examples." (Hilbert)

We will count the number of distinct color patterns of the unit disc in the plane partitioned into four quadrants (1, 2, 3, 4) by the two coordinate axes. Each quadrant is to be colored in one of three colors: **r**ed, **w**hite, or **b**lue. To highlight the different symmetry assumptions one may impose on the disc and the group action appropriate for each situation (and to help understand the meaning of the Theorem), we will do the counting in three different settings.

Setting **I**. The disc is immovable.

In this case, no question of symmetry is involved. Each quadrant can be colored in 3 ways. There are 4 quadrants. The total number of color patterns is $3^4$.

Setting **R**. The disc can be freely rotated about its center.

In this situation, for example, the coloring red in the first quadrant white in the three remaining quadrants, rwww, is indistinguishable from the colorings wrww, wwrw, and wwwr. (These last three colorings can be obtained from the first one by a rotation of $90^o$, $180^o$ and $270^o$, respectively.)  To find the correct number of distinct color patterns, all you need to do is to carry through the following simple procedure.

Take the cyclic permutation group of order 4 generated by the "$90^0$ rotation" (1234). Express each element of the group in its disjoint cycle form.  Do not omit 1-cycles.

C = {(1)(2)(3)(4), (1234), (13)(24), (1432)}.

Change every digit to x and simplify.  Keep all parentheses.  Then take their average (add up the monomials and divide the sum by the order of the group, 4).

$(1/4)((x)^4 + (x^2)^2 + 2(x^4))$.

Substitute $(r^i + w^i + b^i)$ for $x^i$.  You will get

$I_C(r, w, b) = (1/4)((r+w+b)^4 + (r^2+w^2+b^2)^2 + 2(r^4+w^4+b^4))$

$= r^4 + r^3w + br^3 + 2r^2w^2 + 3br^2w + 2b^2r^2 + rw^3 + 3brw^2 + 3b^2rw + b^3r + w^4 + bw^3 + 2b^2w^2 + b^3w + b^4$.

**The coefficient of $r^iw^jb^k$ in $I_C(r, w, b)$ is the number of distinct color patterns with i quadrants colored red, j quadrants colored white, and k quadrants colored blue.**  The total number of distinct color patterns is therefore $= I_C(1, 1, 1) = 24$.

Needless to say, there is much more information contained in $I_C(r, w, b)$.

Setting **RF**.  The disc can be rotated and in addition it can be flipped over and the color is visible on both sides (same as a necklace with 4 beads colored red, white, blue).

The only thing we have to do to accommodate the additional symmetry is to change the permutation group.  We use, in place of the cyclic permutation group of order 4 in Setting R, the dihedral permutation group of order 8 generated by the "$90^0$ rotation" (1234) and the "flip about the x-axis" (14)(23)  Everything else works exactly the same way.

Disjoint cycles:  D = {(1)(2)(3)(4), (1234), (13)(24), (1432), (14)(23), (1)(3)(24),

(12)(34), (13)(2)(4)}.

Monomial average:  $(1/8)((x)^4 + 3(x^2)^2 + 2(x)^2(x^2) + 2(x^4))$.

Substitution $(r^i + w^i + b^i ) \rightarrow x^i$:  $I_D(r, w, b) = (1/8)((r+w+b)^4 + 3(r^2+w^2+b^2)^2 +$

$2(r+w+b)^2(r^2+w^2+b^2) + 2(r^4+w^4+b^4))$.

The total number of distinct color patterns is = $I_D(1, 1, 1) = 21$.

Pick up a book on combinatorics (of some depth). Chances are good you will find in it a chapter devoted to Polya's theory ([1], [2], [3], [4], [6]). An English translation of Polya's original paper is in [5].

I thank Rena Yang for valuable suggestions.